\documentclass[10pt]{amsart}  
\usepackage{amsmath,amssymb,amsfonts,graphicx, hyperref}

\title{High accuracy semidefinite programming bounds for kissing numbers}

\author{Hans D. Mittelmann} 
\address{H.D.~Mittelmann, School of Mathematical
  and Statistical Sciences, Arizona State University, Tempe, AZ
  85287-1804, USA} 
\email{mittelmann@asu.edu}

\author{Frank Vallentin} 
\address{F.~Vallentin, Delft Institute of Applied Mathematics, Technical University of Delft, P.O. Box 5031, 2600 GA Delft, The Netherlands}
\email{f.vallentin@tudelft.nl}

\thanks{The second author was partially supported by
  the Deutsche Forschungsgemeinschaft (DFG) under grant SCHU
  1503/4.}



\date{June 26, 2009}

\newcommand{\defi}[1]{\textit{#1}}

\newcommand{\R}{\mathbb{R}}

\newcommand{\Z}{\mathbb{Z}}
\newcommand{\C}{\mathbb{C}}

\newtheorem{defin}{Definition}[section]

\newtheorem{question}[defin]{Question}

\DeclareMathOperator{\trace}{trace}

\begin{document}

\begin{abstract}
  The kissing number in $n$-dimensional Euclidean space is the maximal
  number of non-overlapping unit spheres which simultaneously can touch
  a central unit sphere. Bachoc and Vallentin developed a method to
  find upper bounds for the kissing number based on semidefinite
  programming. This paper is a report on high accuracy calculations of
  these upper bounds for $n \leq 24$. The bound for $n = 16$ implies a
  conjecture of Conway and Sloane: There is no $16$-dimensional periodic sphere
  packing with average theta series
\[
1 + 7680q^3 + 4320q^4 + 276480q^5 + 61440q^6 + \cdots
\]
\end{abstract}

\maketitle

\markboth{H.D.~Mittelmann, F.~Vallentin}{High accuracy semidefinite
  programming bounds for kissing numbers}

\section{Introduction}
\label{sec:introduction}

In geometry, the \textit{kissing number} in $n$-dimensional Euclidean
space is the maximal number of non-overlapping unit spheres which
simultaneously can touch a central unit sphere. The kissing number is
only known in dimensions $n = 1, 2, 3, 4, 8, 24$, and there were many
attempts to find good lower and upper bounds. We refer to Casselman
\cite{C} for the history of this problem and to Pfender, Ziegler
\cite{PZ}, Elkies \cite{E}, and Conway, Sloane \cite{CS} for more
background information on sphere packing problems.

Bachoc and Vallentin \cite{BV} develop a method
(Section~\ref{sec:notation} recalls it) to find upper bounds for the
kissing number based on semidefinite programming.  Table~1 in
Section~\ref{sec:bounds}, the main contribution of this paper, gives
the values, i.e.\ the first ten significant digits, of these upper
bounds for all dimensions $3, \ldots, 24$. In all cases they are the
best known upper bounds. Dimension $5$ is the first
dimension in which the kissing number is not known. With our
computation we could limit the range of possible values from $\{40,
\ldots, 45\}$ to $\{40, \ldots, 44\}$. In Section~\ref{sec:d4} we show
that the high accuracy computations for the upper bounds in dimension
$4$ result into a question about a possible approach to prove the
uniqueness of the kissing configuration in $4$ dimensions.

Although acquiring the data for the table is a purely computational
task we think that providing this table is valuable for several
reasons: The kissing number is an important constant in geometry and
results can depend on good upper bounds for it. For instance in
Section~\ref{sec:nonexistence} we show that there is no periodic point set
in dimension $16$ with average theta series
\begin{equation*}
1 + 7680q^3 + 4320q^4 + 276480q^5 + 61440q^6 + \cdots
\end{equation*}
This proves a conjecture of Conway and Sloane \cite[Chapter 7, page
190]{CS}.  Furthermore, the actual computation of the table was very
challenging. Bachoc and Vallentin \cite{BV} gave results for dimension
$3, \ldots, 10$. However, they report on numerical difficulties which
prevented them from extending their results. Now using new, more
sophisticated high accuracy software and faster computers and more
computation time we could overcome some of the numerical
difficulties. Section~\ref{sec:bounds} contains details about the
computations.

\section{Notation}
\label{sec:notation}

In this section we set up the notation which is needed for our
computation. For more information we refer to \cite{BV}. For natural
numbers $d$ and $n \geq 3$ let $s_d(n)$ be the optimal value of the
minimization problem
\begin{eqnarray*}
\begin{split}
  \min\Big\{ & 1 + \sum_{k=1}^d a_k  + b_{11} + \langle
  F_0, S_0^n(1,1,1) \rangle : \\
& a_1, \ldots, a_d \in \R,\;\; a_1, \ldots, a_d \geq 0,\\
& b_{11}, b_{12}, b_{22} \in \R,\;\; \left(\begin{smallmatrix} b_{11} & b_{12}\\ b_{12} & b_{22} \end{smallmatrix}\right) \text{ is positive semidefinite, }\\
& F_k \in \R^{(d+1-k) \times (d+1-k)},\;\; F_k \text{ is positive semidefinite, }\;\; k = 0, \ldots, d,\\
& q, q_1 \in \R[u],\;\; \deg (p + pq_1) \leq d,\;\; \text{$p, p_1$ sums of squares},\\
& r, r_1, \ldots, r_4 \in \R[u,v,t],\;\; \deg (r + \sum_{i=1}^4 p_ir_i) \leq d,\;\; \text{$r, r_1, \ldots, r_4$ sums of squares},\\
& 1 + \sum_{k=1}^d a_k P_k^n(u) + 2b_{12} +
b_{22} +3 \sum_{k=0}^d \langle F_k, S_k^n(u,u,1)\rangle + q(u) + p(u)q_1(u)= 0,\\
& b_{22} + \sum_{k=0}^d \langle F_k, S_k^n(u, v,t)\rangle + r(u,v,t) + \sum_{i=1}^4 p_i(u,v,t) r_i(u,v,t)= 0 \big\}.
\end{split}
\end{eqnarray*}
Here $P^n_k$ is the normalized Jacobi polynomial of degree $k$ with
$P^n_k(1) = 1$ and parameters $((n-3)/2,(n-3)/2)$. In general,
\textit{Jacobi polynomials} with parameters $(\alpha, \beta)$ are
orthogonal polynomials for the measure $(1-u)^{\alpha}(1+u)^{\beta}
du$ on the interval $[-1,1]$.  Before we can define the matrices
$S^n_k$ we first define the entry $(i,j)$ with $i, j \geq 0$ of the
(infinite) matrix $Y^n_k$ containing polynomials in the variables
$u,v,w$ by
\begin{equation*}
\begin{split}
\big(Y^n_k\big)_{i,j}(u,v,t) & = u^i v^j \;\cdot\\ &((1-u^2)(1-v^2))^{k/2}
P^{n-1}_k\left(\frac{t-uv}{\sqrt{(1-u^2)(1-v^2)}}\right).
\end{split}
\end{equation*}
Then we get $S^n_k$ by symmetrization: $S^n_k =
\frac{1}{6}\sum_{\sigma} \sigma Y^n_k$, where $\sigma$ runs through
all permutations of the variables $u,v,t$ which acts on the matrix
coefficients in the obvious way.  The polynomials $p$, $p_1, \ldots,
p_4$ are given by
\begin{equation*}
\begin{split}
& p(u) = -(u+1)(u+1/2),\\
& p_1(u,v,t) = p(u),\quad p_2(u,v,t) = p(v),\quad p_3(u,v,t) = p(t),\\
& p_4(u,v,t) = 1 + 2uvt - u^2 - v^2 - t^2.
\end{split}
\end{equation*}
By $\langle A, B \rangle$ we denote the inner product between
symmetric matrices $\trace(AB)$.

In \cite{BV} it is shown that this minimization problem is a
semidefinite program and that every upper bound on $s_d(n)$ provides
an upper bound for the kissing number in dimension $n$. Clearly, the
numbers $s_d(n)$ form a monotonic decreasing sequence in $d$.

\section{Bounds for kissing numbers}
\label{sec:bounds}

\medskip

\begin{center}
\small
\begin{tabular}{c|c|c|c}
    & best lower  & best upper bound &  SDP \\
$n$ & bound known & previously known &  bound \\
\hline
3  &     12 & \underline{12} & $s_{11}(3) = 12.42167009\dots$ \\
   &        & \cite{SW} Sch\"utte, v.d.\ Waerden, 1953 & $s_{12}(3) = 12.40203212\dots$\\
   &        &                                          & $s_{13}(3) = 12.39266509\dots$\\
   &        &                                          & $s_{14}(3) = 12.38180947\dots$\\
\hline
4  &     24 & \underline{24}         & $s_{11}(4) = 24.10550859\dots$\\
   &        & \cite{M} Musin, 2008 & $s_{12}(4) = 24.09098111\dots$\\
   &        &                                          & $s_{13}(4) = 24.07519774\dots$\\
   &        &                                          & $s_{14}(4) = 24.06628391\dots$\\
\hline
5  &     40 & 45   &  $s_{11}(5) = 45.06107293\dots$\\
   &        & \cite{BV} Bachoc, Vallentin, 2008 & $s_{12}(5) = 45.02353644\dots$\\
   &        &                                          & $s_{13}(5) = 45.00650838\dots$\\
   &        &                                          & $s_{14}(5) = \underline{44}.99899685\dots$\\
\hline
6  &     72 & \underline{78}  &  $s_{11}(6) = 78.58344077\dots$\\
   &        & \cite{BV} Bachoc, Vallentin, 2008 & $s_{12}(6) = 78.35518719\dots$\\
   &        &                                          & $s_{13}(6) = 78.29404232\dots$\\
   &        &                                          & $s_{14}(6) = 78.24061272\dots$\\
\hline
7  &    126 & 135& $s_{11}(7) = \underline{134}.8824614\dots$\\
   &        & \cite{BV} Bachoc, Vallentin, 2008 & $s_{12}(7) = 134.7319671\dots$\\
   &        &                                          & $s_{13}(7) = 134.5730609\dots$\\
   &        &                                          & $s_{14}(7) = 134.4488169\dots$\\
\hline
8  &    240 & \underline{240}  & $s_{11}(8) = 240.0000000\dots$ \\
   &        & \cite{OS} Odlyzko, Sloane, 1979 & \\
   &        & \cite{Le} Levenshtein, 1979 & \\
\hline
9  &    306 & 366 & $s_{11}(9) = 365.3229274\dots$ \\ 
   &        & \cite{BV} Bachoc, Vallentin, 2008 & $s_{12}(9) = \underline{364}.7282746\dots$\\
   &        &                                          & $s_{13}(9) = 364.3980087\dots$\\
   &        &                                          & $s_{14}(9) = 364.0919287\dots$\\
\hline
10 &    500 & 567\hspace{1ex}  & $s_{11}(10) = 558.1442813\dots$ \\
   &        & \cite{BV} Bachoc, Vallentin, 2008 & $s_{12}(10) = 556.2840736\dots$\\
   &        &                                          & $s_{13}(10) = 555.2399024\dots$\\
   &        &                                          & $s_{14}(10) = \underline{554}.5075418\dots$\\
\hline
11 &    582 & 915                   & $s_{11}(11) = 878.6158044\dots$ \\
   &        & \cite{OS} Odlyzko, Sloane, 1979 & $s_{12}(11) = 873.3790094\dots$\\
   &        &                                          & $s_{13}(11) = 871.9718533\dots$\\
   &        &                                          & $s_{14}(11) = \underline{870}.8831157\dots$\\
\hline
\end{tabular}
\end{center}

\begin{center}
\small
\begin{tabular}{c|c|c|c}
12 &    840 & 1416                & $s_{11}(12) = 1364.683810\dots$  \\
   &        & \cite{OS} Odlyzko, Sloane, 1979 & $s_{12}(12) = 1362.200297\dots$ \\
   &        &                                          & $s_{13}(12) = 1359.283834\dots$\\
   &        &                                          & $s_{14}(12) = \underline{1357}.889300\dots$\\
\hline
13 &   1130 & 2233                   & $s_{11}(13) = 2089.116331\dots$  \\
   &        & \cite{OS} Odlyzko, Sloane, 1979 & $s_{12}(13) = 2080.631518\dots$\\
   &        &                                          & $s_{13}(13) = 2073.074796\dots$\\
   &        &                                          & $s_{14}(13) = \underline{2069}.587585\dots$\\
\hline
14 &   1582 & 3492                  & $s_{11}(14) = 3224.950751\dots$ \\
   &        & \cite{OS} Odlyzko, Sloane, 1979 & $s_{12}(14) = 3202.448902\dots$\\
   &        &                                          & $s_{13}(14) = 3189.127644\dots$\\
   &        &                                          & $s_{14}(14) = \underline{3183}.133169\dots$\\
\hline
15 &   2564 & 5431                   & $s_{11}(15) = 4949.650431\dots$ \\
   &        & \cite{OS} Odlyzko, Sloane, 1979 & $s_{12}(15) = 4893.479446\dots$\\
   &        &                                          & $s_{13}(15) = 4876.037229\dots$\\
   &        &                                          & $s_{14}(15) = \underline{4866}.245659\dots$\\
\hline
16 &   4320 & 8312                   & $s_{11}(16) = 7515.952644\dots$ \\
   &        & \cite{P} Pfender, 2007 & $s_{12}(16) = 7432.720718\dots$\\
   &        &                                          & $s_{13}(16) = 7374.093742\dots$\\
   &        &                                          & $s_{14}(16) = \underline{7355}.809036\dots$\\
\hline
17 &   5346 & 12210 &  $s_{11}(17) = 11568.41674\dots$\\
   &        & \cite{P} Pfender, 2007 & $s_{12}(17) = 11333.84265\dots$\\
   &        &                                          & $s_{13}(17) = 11128.26227\dots$\\
   &        &                                          & $s_{14}(17) = \underline{11072}.37543\dots$\\
\hline
18 &   7398 & 17877  & $s_{11}(18) = 17473.48016\dots$ \\
   &        & \cite{OS} Odlyzko, Sloane, 1979 & $s_{12}(18) = 17034.32488\dots$\\
   &        &                                          & $s_{13}(18) = 16686.28908\dots$\\
   &        &                                          & $s_{14}(18) = \underline{16572}.26478\dots$\\
\hline
19 &  10668 & 25900  & $s_{11}(19) = 26397.34794\dots$\\
   &        & \cite{Boy} Boyvalenkov, 1994 & $s_{12}(19) = 25636.98958\dots$ \\
   &        &                                          & $s_{13}(19) = 25029.87432\dots$\\
   &        &                                          & $s_{14}(19) = \underline{24812}.30254\dots$\\
\hline
20 &  17400 & 37974  & $s_{11}(20) = 39045.32761\dots$ \\
   &        & \cite{OS} Odlyzko, Sloane, 1979 & $s_{12}(20) = 37844.10380\ldots$\\
   &        &                                          & $s_{13}(20) = 37067.18966\dots$\\
   &        &                                          & $s_{14}(20) = \underline{36764}.40138\dots$\\
\hline
21 &  27720 & 56851  & $s_{11}(21) = 58087.03857\dots$ \\
   &        & \cite{Boy} Boyvalenkov, 1994 & $s_{12}(21) = 56079.21685\ldots$\\
   &        &                                          & $s_{13}(21) = 55170.03449\dots$\\
   &        &                                          & $s_{14}(21) = \underline{54584}.76757\dots$\\
\hline
22 &  49896 & 86537  &  $s_{11}(22) = 87209.06261\dots$\\
   &        & \cite{OS} Odlyzko, Sloane, 1979 & $s_{12}(22) = 84922.09101\dots$\\
   &        &                                          & $s_{13}(22) = 84117.92103\dots$\\
   &        &                                          & $s_{14}(22) = \underline{82340}.08003\dots$\\
\hline
23 &  93150 & 128095 &  $s_{11}(23) = 128360.7969\dots$\\
   &        & \cite{Boy} Boyvalenkov, 1994 & $s_{12}(23) = 127323.7095\dots$\\
   &        &                                          & $s_{13}(23) = 125978.7655\dots$\\
   &        &                                          & $s_{14}(23) = \underline{124416}.9796\dots$\\
\hline
24 & 196560 & \underline{196560} &  $s_{11}(24) = 196560.0000\dots$\\
   &        & \cite{OS} Odlyzko, Sloane, 1979 & \\
   &        & \cite{Le} Levenshtein, 1979 & \\
\hline
\end{tabular}

\medskip

{\sc Table 1.}
New upper bounds for the kissing number (best known values are underlined). 
\end{center}

Finding the solution of the semidefinite program defined in
Section~\ref{sec:notation} is a computational challenge. It turns out
that the major obstacle is numerical instability and not the problem
size. When $d$ is fixed, then the size of the input matrices stays
constant with $n$; when $n$ is fixed, then it grows rather
moderately with $d$.

There is a number of available software packages for solving
semidefinite programs. Mittelmann compares many existing packages in
\cite{Mi}. For our purpose first order, gradient-based methods such as
\texttt{SDPLR} are far too inaccurate, and second order, primal-dual
interior point methods are more suitable. Here increasingly ill-conditioned
linear systems have to be solved even if the underlying problem is
well-conditioned. This happens in the final phase of the algorithm
when one approaches an optimal solution. Our problems are not
well-conditioned and even the most robust solver \texttt{SeDuMi} which
uses partial quadruple arithmetic in the final phase does not produce
reliable results for $d > 10$.

We thus had to fall back on the implementation \texttt{SDPA-GMP}
\cite{SDPA} which is much slower but much more accurate than other
software packages because it uses the GNU Multiple Precision
Arithmetic Library. We worked with $200$--$300$ binary digits and
relative stopping criteria of $10^{-30}$. The ten significant digits
listed in the table are thus guaranteed to be correct. One problem was
the convergence. Even with the control parameter settings recommended
by the authors of \texttt{SDPA-GMP} for ``\texttt{slow but stable}''
computations, the algorithm failed to converge in several
instances. However, we found parameter choices which worked for all
cases: We varied the parameter \texttt{lambdaStar} between $100$ and
$10000$ depending on the case while the other parameters could be
chosen at or near the values recommended for ``\texttt{slow but
  stable}'' performance.

The computations were done on Intel Core 2 platforms with one and two
Quad processors. The computation time varied between five and ten
weeks per case for $d = 12$. An accuracy of 128 bits in
\texttt{SDPA-GMP} did yield sufficient accuracy but did not yield a
reduction in computing time.

After the computations for the cases $d = 11$ and $d = 12$ were
finished new 128-bit versions (quadruple precision) of \texttt{SDPA}
and \texttt{CSDP} became available; partly with our assistance. These
new versions do not rely on the GNU Multiple Precision Arithmetic
Library. So the computation time for the cases $d = 13$ and $d = 14$
were reasonable: from one week to two and a half weeks.

\section{Question about the optimality of the $D_4$ root system}
\label{sec:d4}

Looking at the values $s_d(4)$ in Table~1 one is led to the following question:

\begin{question}
  Is $\lim\limits_{d \to \infty} s_d(4) = 24$?
\end{question}

If the answer to this question is yes (which at the moment appears unlikely because we computed $s_{15}(4) = 23.06274835\dots$), then it would have two noteworthy consequences about optimality properties of the root system $D_4$.

The root system $D_4$ defines (up to orthogonal transformations) a
point configuration on the unit sphere $S^3 = \{x \in \R^4 : x \cdot x
= 1\}$ consisting of $24$ points; it is the same point configuration
as the one coming from the vertices of the regular $24$ cell. This
point configuration has the property that the spherical distance of
every two distinct points is at least $\arccos 1/2$. Hence, these points
can be the maximal $24$ touching points of unit spheres kissing the
central unit sphere $S^3$.

If $\lim_{d \to \infty} s_d(4) = 24$, then this would prove that the
root system $D_4$ is the unique optimal point configuration of
cardinality~$24$. Here optimality means that one cannot distribute
$24$ points on $S^3$ so that the minimal spherical distance between
two distinct points exceeds~$\arccos 1/2$. Thus, the root system $D_4$
would be characterized by its kissing property. This is generally
believed to be true but so far no proof could be given.

Another consequence would be that there is no universally optimal
point configuration of $24$ points in $S^3$ as conjectured by Cohn,
Conway, Elkies, Kumar \cite{CCEK}. Universally optimal point
configurations minimize every absolutely monotonic potential function.
The conjecture would follow if the answer to our question is yes:
Every universally optimal point configuration is automatically optimal
and Cohn, Conway, Elkies, Kumar \cite{CCEK} show that the root system
$D_4$ is not universally optimal.

\section{Nonexistence of a sphere packing}
\label{sec:nonexistence}

Our new upper bound of $7355$ for the kissing number in dimension $16$
implies that there is no periodic point set in dimension $16$ whose
average theta series equals
\begin{equation}
\label{eq:extremal}
1 + 7680q^3 + 4320q^4 + 276480q^5 + 61440q^6 + \cdots.
\end{equation}
This settles a conjecture of Conway and Sloane \cite[Chapter 7, page
190]{CS}. In this section we explain this result. We refer to Conway, Sloane \cite{CS}, Elkies \cite{E}, and to Bowert
\cite{Bow} for more information.

An $n$-dimensional \defi{periodic point set} $\Lambda$ is a finite
union of translates of an $n$-dimensional lattice, i.e.\ one can write
$\Lambda$ as $\Lambda = (A\Z^n + v_1) \cup \ldots \cup (A\Z^n +
v_N)$, with $v_1, \ldots, v_N \in \R^n$, and $A : \R^n \to \R^n$ is a
linear isomorphism. The \defi{average theta series} of $\Lambda$ is
\begin{equation*}
  \Theta_{\Lambda}(z) = \frac{1}{N} \sum_{i=1}^N \sum_{j = 1}^N \sum_{v \in \Z^n} q^{\|Av - v_i + v_j\|^2}, \text{ with } q = e^{\pi i z}.
\end{equation*}
This is a holomorphic function defined on the complex upper
half-plane. A holomorphic function $f$ which is defined on the complex
upper half-plane, which is meromorphic for $z \to i\infty$, and which
satisfies the transformation laws
\begin{equation*}
f(-1/z) = z^8 f(z), \text{ and } f(z + 2) = f(z) \text { for all $z \in \C$ with $\Im z > 0$},
\end{equation*}
is called a \defi{modular forms} of weight $8$ for the Hecke group
$G(2)$. The expression \eqref{eq:extremal} defines the unique modular
form of weight $8$ for the Hecke group $G(2)$ which starts off with $1
+ 0q^1 + 0q^2$. It is also called an \defi{extremal modular form}, see
Scharlau and Schulze-Pillot \cite{SSP}

If there would be a $16$-dimensional periodic point set whose
average theta series coincides with \eqref{eq:extremal} then this
periodic point set would define the sphere centers of a sphere packing
with extraordinary high density (see \cite[Chapter 7, page 190]{CS}). At
the same time the existence of such a periodic point set would show
that the kissing number in dimension $16$ is at least $7680$. This is
not the case.

\section*{Acknowledgements}

We thank Etienne de Klerk and Renata Sotirov for initiating our
collaboration and we thank Frank Bowert and Rudolf Scharlau for
bringing the conjecture of Conway and Sloane to our attention.


\begin{thebibliography}{[99]}

\bibitem{BV} C.~Bachoc, F.~Vallentin, {\em New upper bounds for
    kissing numbers from semidefinite programming},
  J. Amer. Math. Soc. \textbf{21} (2008), 909--924.

\bibitem{Bow} F.~Bowert, {\em Gewichtsz\"ahler und Distanzz\"ahler von
    Codes und Kugelpackungen}, Ph.D. thesis, University of Dortmund,
  2004.

\bibitem{Boy} P.~Boyvalenkov, {\em Small improvements of the upper
    bounds of the kissing numbers in dimensions $19$, $21$ and $23$},
  Atti Sem. Mat. Fis. Univ. Modena {\bf 42} (1994), 159--163.

\bibitem{C} B.~Casselman,
{\em The difficulties of kissing in three dimensions},
Notices Amer. Math. Soc. {\bf 51} (2004), 884--885. 

\bibitem{CCEK} H.~Cohn, J.H.~Conway, N.D.~Elkies, A.~Kumar, {\em The
    $D_4$ root system is not universally optimal}, Exp. Math. {\bf
    16} (2007), 313--320.

\bibitem{CS} J.H.~Conway, N.J.A.~Sloane, {\em Sphere packings, lattices
and groups}, third edition, Springer, 1999.

\bibitem{E} N.D.~Elkies, {\em Lattices, linear codes, and invariants,
  Part I/II}, Notices Amer. Math. Soc. {\bf 47} (2000), 1238--1245,
1382--1391.

\bibitem{SDPA} K.~Fujisawa, M.~Fukuda, K.~Kobayashi, M.~Kojima,
  K.~Nakata, M.~Nakata, M.~Yamashita, {\em SDPA (SemiDefinite
    Programming Algorithm) and SDPA-GMP user's manual --- version
    7.1.1}, Research Report B-448, Department of Mathematical and
  Computing Sciences, Tokyo Institute of Technology, June 2008.

\bibitem{Le} V.I.~Levenshtein, {\em On bounds for packing in
    $n$-dimensional Euclidean space}, Soviet Math. Dokl. {\bf 20}
  (1979), 417--421.

\bibitem{Mi} H.D.~Mittelmann, {\em An independent benchmarking of SDP 
and SOCP solvers},
  Math. Prog. \textbf{95} (2003), 407--430.

\bibitem{M} O.R.~Musin, {\em The kissing number in four
dimensions}, Ann. of Math. \textbf{168} (2008), 1--32.

\bibitem{OS} A.M.~Odlyzko, N.J.A.~Sloane,
{\em New bounds on the number of unit spheres that can touch a unit
  sphere in n dimensions}, 
J. Combin. Theory Ser. A {\bf 26} (1979), 210--214.

\bibitem{P} F.~Pfender, {\em Improved Delsarte bounds for spherical
    codes in small dimensions}, J. Combin. Theory Ser. A {\bf 114}
  (2007), 1133--1147.

\bibitem{PZ} F.~Pfender, G.M.~Ziegler,
{\em Kissing numbers, sphere packings and some unexpected proofs},
Notices Amer. Math. Soc.  {\bf 51}  (2004), 873--883. 

\bibitem{SSP} R.~Scharlau, R.~Schulze-Pillot, {\em Extremal lattices},
  p. 139--170 in {\sl Algorithmic algebra and number theory
    (Heidelberg, 1997)}, Springer, 1999.

\bibitem{SW} K.~Sch\"utte, B.L.~van der Waerden,
{\em Das Problem der dreizehn Kugeln}, 
Math. Ann. {\bf 125} (1953) 325--334.

\end{thebibliography}
\end{document}